\documentclass[11pt,twoside]{amsart}

\usepackage{pstricks}
\usepackage{multido}

\usepackage{graphicx,url}
\usepackage{amssymb}
\usepackage{amsmath}
\usepackage{amscd}
\usepackage{amstext}
\usepackage{amsbsy}
\usepackage{amsxtra}
\usepackage{amsfonts}
\usepackage{latexsym}
\usepackage{alltt}
\usepackage{bbm}

\usepackage{verbatim}
\usepackage{amscd}
\usepackage[all]{xy}

\theoremstyle{plain}
\newtheorem{theorem}{Theorem}
\newtheorem*{theorem*}{Theorem}
\newtheorem{lemma}{Lemma}[section]

\newtheorem{remark}[lemma]{Remark}

\numberwithin{equation}{section}

\title[On the elliptic sinh-Gordon equation...]
{On the elliptic sinh-Gordon equation with Durham boundary conditions}

\author{M. Kilian}
\address{M. Kilian, University College Cork, Ireland.}
\email{m.kilian@ucc.ie}

\author{G. Smith}
\address{G. Smith, Universidade Federal do Rio de Janeiro, Brazil.}
\email{grahamandrewsmith@gmail.com}

\newcommand{\opSinh}{{\text{\rm sinh}}}
\newcommand{\opCosh}{{\text{\rm cosh}}}
\newcommand{\opsl}{{\text{\rm sl}}}
\newcommand{\opDet}{{\text{\rm Det}}}
\newcommand{\opTr}{{\text{\rm Tr}}}
\newcommand{\opAdj}{{\text{\rm Adj}}}
\newcommand{\opId}{{\text{\rm Id}}}
\newcommand{\opIm}{{\text{\rm Im}}}
\newcommand{\opRe}{{\text{\rm Re}}}
\newcommand{\opDim}{{\text{\rm Dim}}}
\def\eqnum#1{#1}
\def\nexteqnno[#1]{\label{#1}}
\def\eqnref#1{\eqref{#1}}
\def\procref#1{\ref{#1}}
\def\subheadref#1{\ref{#1}}
\def\proclabel#1{\label{#1}}
\def\newhead#1[#2]{\begin{center}\sc #1\end{center}}
\def\newsubhead#1[#2]{\section{#1}\label{#2}}
\def\Cal#1{\mathcal{#1}}
\def\myitem#1{\leavevmode\newline\noindent\hbox to .5cm{\hfill#1\hss}}
\catcode`\@=11
\def\eqalign#1{\null\,\vcenter{\openup1\jot \m@th %
\ialign{\strut\hfil$\displaystyle{##}$&$\displaystyle{{}##}$\hfil\crcr#1\crcr}}\,}
\def\eqalignleft#1{\null\,\vcenter{\openup1\jot \m@th %
\ialign{\strut$\displaystyle{##}$\hfil&$\displaystyle{{}##}$\hfil\crcr#1\crcr}}\,}
\def\triplealign#1{\null\,\vcenter{\openup1\jot \m@th %
\ialign{\strut\hfil$\displaystyle{##}\quad$&\hfil$\displaystyle{{}##}$&$\displaystyle{{}##}$\hfil\crcr#1\crcr}}\,}
\catcode`\@=12
\begin{document}


%=========================

\begin{abstract}
We adapt Sklyanin's $K$-matrix formalism to the sinh-Gordon equation, and prove that all free boundary constant mean curvature (CMC) annuli in the unit ball in $\mathbb{R}^3$ are of finite type.
\end{abstract}
\dedicatory{Dedicated to Ulrich Pinkall on the occasion of his 65th birthday.}
\maketitle
%
%\newhead{Introduction}[Introduction]
%
\newsubhead{CMC immersions and the sinh-Gordon equation}[TheSinhGordonEquation]
Let $\Omega\subseteq\Bbb{C}$ be an open subset with smooth boundary $\partial\Omega$. For $r>0$, let $B_r(0)$ denote the open ball of radius $r$ about the origin in $\Bbb{R}^3$. For all $r$, a minimal or CMC immersion $f:\overline{\Omega}\rightarrow B_r(0)$ is said to have {\sl free boundary} whenever it meets $S_r(0):=\partial B_r(0)$ orthogonally along $\partial\Omega$. Free boundary minimal and CMC surfaces have attracted the interest of geometric analysts since the work \cite{FraserSchoenI, FraserSchoenII} of Fraser--Schoen. The purpose of this paper is to open the way to applying integrable systems techniques to the study of free boundary CMC annuli. The free boundary minimal case will not be treated here.
\par
In order to better explain the ideas studied in the sequel, we first review the case without boundary, known as the {\sl bulk} case. Here, the modern use of integrable systems techniques for the study of CMC tori in $\Bbb{R}^3$ traces its roots to the original construction \cite{Wente} of Wente, which was further developed by Abresch in \cite{Abresch}. These ideas were extended to their full generality by Pinkall--Sterling in \cite{PinkallSterling} by showing that, modulo closing conditions, the study of immersed CMC tori in $\Bbb{R}^3$ is equivalent to the study of real, doubly-periodic solutions of the elliptic sinh-Gordon equation. Independently, in \cite{Hitchin}, a similar technique was developed by Hitchin for classifying all harmonic tori in the $3$-sphere. In this paper, we will follow Pinkall--Sterling's approach.
\par
Pinkall--Sterling proceed as follows. Let $f:\Bbb{C}\rightarrow\Bbb{R}^3$ be a smooth, doubly-periodic immersion of non-zero constant mean curvature $H$ (where here we define the mean curvature to be equal to the algebraic mean of the two principal curvatures). We suppose that $f$ is conformal, so that the metric it induces over $\Omega$ is given by
\begin{equation}
g := e^{2\omega}dz d\overline{z},\eqnum{\nexteqnno[ConformalMetric]}
\end{equation}
for some smooth, real-valued function $\omega$, which we call the {\sl conformal factor} of $f$. Recall from \cite{Hopf} that the {\sl Hopf differential}
\begin{equation*}
Q := \phi\, dz dz
\end{equation*}
of $f$ is constant. Thus, upon rescaling the domain and the codomain if necessary, we may suppose that
\begin{equation}
H = \frac{1}{2}\ \text{and}\ \left|Q\right| = \frac{1}{4},\eqnum{\nexteqnno[NormalisingConditions]}
\end{equation}
and the Gauss--Codazzi equations for $f$ are then equivalent to
\begin{equation}
\omega_{z\overline{z}} + \frac{1}{8}\opSinh(2\omega) = 0,\eqnum{\nexteqnno[sinhGordon]}
\end{equation}
which is the {\sl elliptic sinh--Gordon equation}.
\par
Conversely, by the fundamental theorem of surface theory, given $H$ and $Q$ satisfying \eqref{NormalisingConditions} and a doubly-periodic function $\omega:\Bbb{C}\rightarrow\Bbb{R}$ satisfying \eqref{sinhGordon}, there exists, up to rigid motions of $\Bbb{R}^3$, a unique CMC-$(1/2)$ immersion $f:\Bbb{C}\rightarrow\Bbb{R}$ with Hopf differential $Q$ and conformal factor $\omega$. This immersion is, furthermore, a torus provided that two further closing conditions are satisfied. In this manner, we obtain the desired equivalence, modulo closing conditions, between CMC immersed tori on the one hand, and doubly periodic solutions of the elliptic sinh--Gordon equation, on the other.
\par
Significantly, the sinh--Gordon equation is one of the best studied equations of integrable systems theory. So it is that Pinkall--Sterling use the results of this highly developed field to classify all CMC immersed tori in $\Bbb{R}^3$. Subsequently, in \cite{Bobenko}, Bobenko uses this classifiction to derive explicit formulae for all constant mean curvature tori in $\Bbb{R}^3$ in terms of theta functions. With these results in mind, it is natural to ask whether the same techniques can be applied in the free-boundary case.
\newsubhead{Finite type solutions}[FiniteTypeSolutions]
In \cite{Bobenko}, Bobenko observes that the key steps in both Pinkall--Sterling's and Hitchin's work lie in showing that all doubly-periodic solutions of the elliptic sinh--Gordon equation over $\Bbb{C}$ are of finite type. Heuristically, this means that all such solutions are completely determined by polynomial data via a certain ansatz (c.f. \cite{Hitchin}). However, the formal statement of the finite-type property is rather technical, with different authors using different definitions. In this paper, we adopt the perspective of Adler--Kostant--Symes theory (c.f. \cite{BurstallFerusPeditPinkall} and \cite{FordyWood}). Recall first that the elliptic sinh--Gordon equation translates into the integrability condition of the Lax pair
\begin{equation}
\eqalign{
\alpha_z(\lambda,\gamma) &= \frac{1}{2}\omega_z\sigma_0 + \frac{i}{4\lambda}e^\omega\sigma_+ + \frac{i\gamma}{4}e^{-\omega}\sigma_-\ \text{and}\cr
\alpha_{\overline{z}}(\lambda,\gamma) &= -\frac{1}{2}\omega_{\overline{z}}\sigma_0 + \frac{i}{4\gamma}e^{-\omega}\sigma_+ + \frac{i\lambda}{4}e^{\omega}\sigma_-,\cr}\eqnum{\nexteqnno[LaxPair]}
\end{equation}
where
\begin{equation}
\sigma_0 := \begin{pmatrix} 1 & 0\cr 0& -1\cr\end{pmatrix},\
\sigma_+ := \begin{pmatrix} 0 & 1\cr 0& 0\cr\end{pmatrix}\ \text{and}\
\sigma_- := \begin{pmatrix} 0 & 0\cr 1& 0\cr\end{pmatrix},\eqnum{\nexteqnno[SigmaMatricesI]}
\end{equation}
and $\lambda,\gamma\in\Bbb{C}^*$ are non-zero complex parameters which we respectively call the {\sl spectral} and {\sl torsion parameters}.\footnote{The Lax pair already appears implicitly in the work \cite{PinkallSterling} of Pinkall--Sterling. Indeed, whilst the torsion parameter is explicitely mentioned, the spectral parameter also appears implicitely as the independent variable of the generating functions of the sequences studied.} For the purposes of this introduction, a {\sl Killing field}\footnote{A different definition will be used in the sequel (c.f. \eqnref{KillingFieldEquation}, below).} of $\omega$ is defined to be a function $\Phi:\Bbb{C}^*\times\Bbb{C}^*\rightarrow C^\infty(\Omega,\opsl_2(\Bbb{C}))$ which solves the system of partial differential equations
\begin{equation}
d\Phi(\lambda,\gamma) = [\Phi(\lambda,\gamma),\alpha(\lambda,\gamma)],\eqnum{\nexteqnno[KillingFieldEqnIntro]}
\end{equation}
where here
\begin{equation*}
\alpha := \alpha_z dz + \alpha_{\overline{z}}d\overline{z}.
\end{equation*}
We say that a Killing field is {\sl polynomial} whenever it takes the form
\begin{equation}
\Phi(\lambda) = \sum_{(m,n)\in A}\Phi_{m,n}\lambda^m\gamma^n,\eqnum{\nexteqnno[PolynomialKillingFieldIntro]}
\end{equation}
for some {\sl finite} subset $A$ of $\Bbb{Z}\times\Bbb{Z}$ where, for all $(m,n)$, $\Phi_{m,n}:\Omega\rightarrow\opsl_2(\Bbb{C})$ is a smooth function. According to Adler--Kostant--Symes theory, a solution $\omega$ of the elliptic sinh--Gordon equation is of {\sl finite type} whenever it admits a polynomial Killing field.
\newsubhead{Free boundary CMC annuli}[FreeBoundaryCMC]
With Bobenko's observation in mind, we wish to show that the conformal factors of free boundary CMC annuli are also finite type solutions of \eqnref{sinhGordon}. Consider therefore a periodic, conformal CMC immersion $f$ defined over the ribbon $\Omega:=\Bbb{R}\times[-T,T]$ with free boundary in $B_r(0)$. By classical surface theory, the Hopf differential $Q$ of $f$ has constant argument along each of the boundary components. It follows by the Schwarz reflection principle that it extends to a bounded, holomorphic form over the whole of $\Bbb{C}$ which, by Liouville's theorem, is constant. By classical surface theory again, the conformal factor $\omega$ satisfies the non-linear boundary condition
\begin{equation}
\omega_y = \frac{\epsilon}{r}e^\omega,\eqnum{\nexteqnno[FreeBoundaryCondition]}
\end{equation}
where $\epsilon$ is equal to $+1$ along the upper boundary component and $-1$ along the lower boundary component. More generally, non-linear boundary conditions of the form
\begin{equation}
\omega_y = Ae^\omega + Be^{-\omega},\eqnum{\nexteqnno[DurhamBoundaryConditions]}
\end{equation}
where $A$ and $B$ are constant along each boundary component, are known as {\sl Durham boundary conditions} (see \cite{CorriganEtAlA, CorriganEtAlB}).
\par
We show
\begin{theorem}\label{IntroFiniteType}
\noindent If $\omega:\Bbb{R}\times[-T,T]\rightarrow\Bbb{R}$ is a singly-periodic solution of the elliptic sinh--Gordon equation with Durham boundary conditions, then $\omega$ is of finite type.
\end{theorem}
Our proof also yields deeper information about the structure of the polynomial Killing fields of such $\omega$. First, we say that a Killing field $\Phi$ satisfies the {\sl Sklyanin condition for fields} whenever, at every point of $\partial\Omega$,
\begin{equation}
K(\lambda,\gamma)\Phi(\lambda,\gamma) = \overline{\Phi(\overline{\lambda},\overline{\gamma})}^t K(\lambda,\gamma),\eqnum{\nexteqnno[SklyaninConditionIntro]}
\end{equation}
for all $\lambda,\gamma\in\Bbb{C}^*$, where
\begin{equation}
K(\lambda,\gamma) := \begin{pmatrix} 4A\gamma-4B\lambda&{\displaystyle{\frac{\lambda}{\gamma} - \frac{\gamma}{\lambda}}}\cr
{\displaystyle{\frac{\lambda}{\gamma} - \frac{\gamma}{\lambda}}}&{\displaystyle{\frac{4A}{\gamma} - \frac{4B}{\lambda}}}\cr\end{pmatrix}\eqnum{\nexteqnno[KMatrixIntro]}
\end{equation}
is Sklyanin's $K$-matrix (see \cite{SklyaninI, SklyaninII}, see also \cite{MacIntyre} for an excellent treatment of Sklyanin's ideas). Theorem \procref{IntroFiniteType} immediately follows from
\begin{theorem}\label{MainTheorem}
\noindent If $\omega:\Bbb{R}\times[-T,T]\rightarrow\Bbb{R}$ is a periodic solution of the elliptic sinh--Gordon equation, then $\omega$ satisfies the Durham boundary conditions if and only if it admits a polynomial Killing field $\Phi$ which satisfies the Sklyanin condition for fields.
\end{theorem}
{\bf\noindent Remark:\ }Theorem \procref{MainTheorem} is proven in Section \subheadref{DesChampsDeKillingPolynomes}, below.
%
%\newhead{The Sklyanin condition}[TheSklyaninCondition]
%
\newsubhead{Killing fields}[LesChampsDeKilling]
As in the introduction, let $\Omega:=\Bbb{R}\times[-T,T]$ and let $\omega:\Omega\rightarrow\Bbb{R}$ be a real solution of the elliptic sinh--Gordon equation. Throughout the current chapter, we will consider $\lambda$ as a variable and $\gamma$ as a constant. We begin by describing the algebraic formalism that we will use in the sequel. Let $X$ be a manifold. Let $E$ be a complex vector space. For $k\in\Bbb{Z}$, a {\sl Laurent series} of degree $k$ over $X$ taking values in $E$ is defined to be a {\sl formal} series of the form
\begin{equation*}
\Phi(\lambda) := \sum_{m=k}^\infty\Phi_m\lambda^m,
\end{equation*}
where, for all $m$, $\Phi_m$ is a smooth function over $X$ taking values in $E$ and $\Phi_k$ is non-zero. We say that the trivial series $\Phi=0$
is a {\sl Laurent series} of degree $+\infty$. Let $\Cal{L}(X,E)$ denote the space of Laurent series over $X$ taking values in $E$. In addition, denote
\begin{equation}
\Cal{L}(X) := \Cal{L}(X,\Bbb{C})\ \text{and}\ \Cal{L} := \Cal{L}(\left\{x\right\},\Bbb{C}),\eqnum{\nexteqnno[SpecialSpacesOfLaurentSeries]}
\end{equation}
where $\left\{x\right\}$ here denotes the manifold consisting of a single point. We readily verify
\begin{lemma}\label{BasicAlgebraicPropertiesOfL}
\leavevmode\null\medskip
\myitem{(1)} $\Cal{L}(X,E)$ is a complex vector space;
\medskip
\myitem{(2)} $\Cal{L}(X)$ is a unitary commutative algebra;
\medskip
\myitem{(3)} $\Cal{L}$ is an algebraic field; and
\medskip
\myitem{(4)} $\Cal{L}(X,E)$ is a module over $\Cal{L}(X)$ and a vector space over $\Cal{L}$.
\end{lemma}
\par
For $k\leq l\in\Bbb{N}$, a {\sl Laurent polynomial} of bidegree $(k,l)$ over $X$ is defined to be a Laurent series of the form
\begin{equation*}
\Phi(\lambda) := \sum_{m=k}^l\Phi_m\lambda^m,
\end{equation*}
where $\Phi_k$ and $\Phi_l$ are non-zero. As before, we call the trivial series $\Phi=0$ a {\sl Laurent polynomial} of bidegree $(+\infty,+\infty)$. Let $\Cal{P}(X,E)$ denote the space of Laurent polynomials over $X$ taking values in $E$.
\par
Observe now that the Lax pair \eqnref{LaxPair} is a $1$-form over $\Omega$ which maps vector fields into $\Cal{P}(\Omega,\opsl(2,\Bbb{C}))$. We define a {\sl Killing field}\footnote{We underline that this definition, which will be used throughout the sequel, is different from that given in the introduction. In \cite{BurstallFerusPeditPinkall}, the authors call such objects {\sl formal Killing fields}. In the present context, the adjective, ``formal'' is superfluous, and we therefore omit it.} of $\omega$ over $\Omega$ to be a Laurent series $\Phi\in\Cal{L}(\Omega,\opsl(2,\Bbb{C}))$ which verifies
\begin{equation}
d\Phi = [\Phi,\alpha]\eqnum{\nexteqnno[KillingFieldEquation]}
\end{equation}
at every point of $\Omega$ (c.f. \cite{BurstallFerusPeditPinkall}). Upon evaluating \eqnref{KillingFieldEquation} term by term, we obtain (see \cite{HauswirthKillianSchmidt})
\begin{lemma}\label{LaxFormulae}
\noindent Let $\omega$ be a solution of the sinh--Gordon equation. Let
\begin{equation*}
\Phi := \sum_{m=k}^\infty\begin{pmatrix} u_m& e^\omega t_m\cr e^\omega s_m& -u_m\cr\end{pmatrix} \lambda^m
\end{equation*}
be a Laurent series over $\Omega$ taking values in $\opsl(2,\Bbb{C})$. Then $\Phi$ is a Killing field if and only if, at every point of $\Omega$ and for all $m$,
\begin{align}
4u_{m,z} + i e^{2\omega} s_{m+1} - i\gamma t_m &= 0,\nexteqnno[RI]\\
4u_{m,\overline{z}} + i\gamma^{-1}s_m - i e^{2\omega}t_{m-1} &= 0,\nexteqnno[RII]\\
4\omega_z t_m + 2 t_{m,z} - i u_{m+1} &= 0,\nexteqnno[RIII]\\
2e^\omega t_{m,\overline{z}} - i \gamma^{-1}e^{-\omega} u_m &= 0,\nexteqnno[RIV]\\
2e^\omega s_{m,z} + i \gamma e^{-\omega} u_m &= 0\ \text{and}\nexteqnno[RV]\\
4\omega_{\overline{z}} s_m + 2 s_{m,\overline{z}} + i u_{m-1} &= 0.\nexteqnno[RVI]
\end{align}
\end{lemma}
\newsubhead{The space of Killing fields}[LEspaceDeChampsDeKilling]
Let $\Cal{K}(\Omega)$ denote the space of Killing fields of $\omega$ over $\Omega$. By \cite{PinkallSterling}, this space is non-trivial. Indeed, in that paper, Pinkall--Sterling construct an explict, non-trivial Killing field of $\omega$ over $\Omega$, which we henceforth refer to as the {\sl Pinkall--Sterling field}\footnote{We remark that the formalism of \cite{PinkallSterling} is slightly different from our own, but can be transformed into our own upon replacing their variable $z$ with the variable $\zeta:=iz/2$ so that
\begin{equation*}
\partial_\zeta = -2i\partial_z\ \text{and}\ \partial_{\overline{\zeta}} = 2i\partial_{\overline{z}}.
\end{equation*}}.
It will be useful to recall the recursive formula used to construct this field. First, define
\begin{equation}
u_0 := 0\ \text{and}\ \psi_0 := -\frac{1}{2}.\eqnum{\nexteqnno[InitialValuesOfRecursion]}
\end{equation}
Next, after having determined $u_1,\cdots,u_m$ and $\psi_1,\cdots,\psi_{m-1}$, define
\begin{equation}
\psi_m :=
\left\{
\eqalignleft{
\gamma u_k^2 + 2\sum_{n=1}^{k-1}\theta_{n,m-n}&\qquad\text{if}\ m=2k-1,\ \text{and}\cr
\gamma u_k u_{k+1} + \theta_{k,k} + 2\sum_{n=1}^{k-1}\theta_{n,m-n}&\qquad\text{if}\ m=2k,\cr}\right.
\eqnum{\nexteqnno[RecursionForPsiI]}
\end{equation}
and
\begin{equation}
u_{m+1} := \frac{1}{\gamma}\big(-4u_{m,zz} + 4i\omega_z\psi_m\big),\eqnum{\nexteqnno[RecursionForU]}
\end{equation}
where, for all $p$ and for all $q$,
\begin{equation}
\theta_{p,q} := \gamma u_p u_{q+1} + 4 u_{p,z}u_{q,z} + \psi_p\psi_q.\eqnum{\nexteqnno[RecursionForPsiII]}
\end{equation}
The sequences $(t_m)$ and $(s_m)$ are now determined by
\begin{equation}
t_m := \frac{1}{\gamma}\big(-2iu_{m,z} - \psi_m)\ \text{and}\
s_m := e^{-2\omega}\big(2iu_{m-1,z} - \psi_{m-1}\big).\eqnum{\nexteqnno[RecursionForTAndS]}%
\end{equation}
The Pinkall--Sterling field is then the series
\begin{equation}
\Phi := \sum_{m=0}^\infty\begin{pmatrix} u_m& e^\omega t_m\cr e^\omega s_m& -u_m\cr\end{pmatrix}\lambda^m.\eqnum{\nexteqnno[ExpansionForPhi]}
\end{equation}
\par
We now show that $\Cal{K}(\Omega)$ is $1$-dimensional over $\Cal{L}$. To this end, we show
\begin{lemma}\proclabel{VanishingConditions}
\noindent Let
\begin{equation*}
\Phi := \sum_{m=0}^\infty\begin{pmatrix} u_m& e^\omega t_m\cr e^\omega s_m& -u_m\cr\end{pmatrix}\lambda^m
\end{equation*}
be a Killing field of $\omega$ over $\Omega$. For all $m$,
\medskip
\myitem{(1)} if $u_m=0$, then $t_m$ and $s_m$ are constant;
\medskip
\myitem{(2)} if $u_m=t_m=0$, then $s_{m+1}=u_{m+1}=0$; and
\medskip
\myitem{(3)} if $u_m=s_m=0$, then $u_{m-1}=t_{m-1}=0$.
\end{lemma}
\proof Suppose that $u_m=0$. By \eqnref{RIV},
\begin{equation*}
t_{m,\overline{z}} = 0.
\end{equation*}
Next, by \eqnref{RI},
\begin{equation*}
e^\omega s_{m+1} = \gamma e^{-\omega}t_m,
\end{equation*}
by \eqnref{RV},
\begin{equation*}
\eqalign{
\gamma e^{-\omega}u_{m+1} &= 2ie^\omega s_{m+1,z}=2i(e^\omega s_{m+1})_z - 2i\omega_z e^\omega s_{m+1}\cr
&=2i\gamma(e^{-\omega}t_m)_z - 2i\gamma\omega_z e^{-\omega}t_m=2i\gamma e^{-\omega}t_{m,z} - 4i\gamma\omega_z e^{-\omega}t_m,\cr}
\end{equation*}
and by \eqnref{RIII},
\begin{equation*}
t_{m,z} = 0.
\end{equation*}
It follows that $t_m$ is constant. In the same manner, we show that $s_m$ is constant, and $(1)$ follows. If $u_m=t_m=0$, then it follows by \eqnref{RI} and \eqnref{RIII} that $s_{m+1}=u_{m+1}=0$, and $(2)$ follows. If $s_m=u_m=0$, then it follows by \eqnref{RII} and \eqnref{RVI} that $u_{m-1}=t_{m-1}=0$, and $(3)$ follows. This completes the proof.\qed
\begin{lemma}\proclabel{Unidimensional}
\noindent $\Cal{K}(\Omega)$ is the $1$-dimensional vector space over $\Cal{L}$ generated by the Pinkall--Sterling field.
\end{lemma}
\remark We see, in particular, that in the present framework a solution of the elliptic sinh--Gordon equation is of finite type if and only if its space of Killing fields is generated by a Laurent polynomial.
\medskip
\proof Let $\Phi$ be an arbitrary Killing field of $\omega$ over $\Omega$ and let $\Psi$ be the Pinkall--Sterling field. We construct recursively a Laurent series $f$ such that
\begin{equation*}
\Phi = f\Psi.
\end{equation*}
Let $k$ be the degree of $\Phi$. As $\Psi$ has degree $0$, the series $f$ must also be of degree $k$. Suppose now that we have already determined the coefficients $f_k,f_{k+1},...,f_{k+l-1}$ in such a manner that the series
\begin{equation*}
\tilde{\Phi} := \Phi - f_{(l)}\Psi
\end{equation*}
is of degree $k+l$, where
\begin{equation*}
f_{(l)} := \sum_{m=k}^{k+l-1}f_m\lambda^m.
\end{equation*}
For all $m$, denote
\begin{equation*}
\tilde{\Phi}_m := \begin{pmatrix} u_m & e^\omega\tau_m \cr e^\omega\sigma_m & -u_m \cr\end{pmatrix}.
\end{equation*}
As $\tilde{\Phi}$ is also a solution of the Killing field equation, it follows by Lemma \procref{VanishingConditions} that
\begin{equation*}
s_{k+l}=0,\ u_{k+l} =0\ \text{and}\ \tau_{k+l} = c,
\end{equation*}
where $c$ is a constant. The result now follows upon setting $f_{k+l}:=2c$.\qed
\newsubhead{The determinant}[LeDeterminant]
We now characterise the Pinkall--Sterling field amongst all Killing fields of $\omega$ over $\Omega$. Observe first that, since elements of $\Cal{K}(\Omega)$ are $2\times 2$ matrices with coefficients in $\Cal{L}(\Omega)$, they have well-defined determinants which are also elements of $\Cal{L}(\Omega)$.
\begin{lemma}
\noindent For every Killing field $\Phi$ of $\omega$ over $\Omega$, $\opDet(\Phi)$ is constant over $\Omega$, that is,
\begin{equation*}
\opDet(\Phi)\in\Cal{L}.
\end{equation*}
\end{lemma}
\proof Indeed
\begin{equation*}
d\opDet(\Phi) = \opTr(\opAdj(\Phi)d\Phi) = \opTr(\opAdj(\Phi)[\Phi,\alpha]) = 0,
\end{equation*}
as desired.\qed
\begin{lemma}%\proclabel{PSGuage}
\noindent For all $\gamma$, the Pinkall--Sterling field is, up to a choice of sign, the unique Killing field $\Phi$ of $\omega$ over $\Omega$ such that
\begin{equation}
\opDet(\Phi) = -\frac{\lambda}{4\gamma}.\eqnum{\nexteqnno[PSGuage]}
\end{equation}
\end{lemma}
\proof Let $\Phi$ be the Pinkall--Sterling field. We first show uniqueness. Let $\Psi$ be another Killing field of $\omega$ over $\Omega$ which satisfies \eqnref{PSGuage}. Since $\Cal{K}(\Omega)$ is generated by $\Phi$, there exists a Laurent series $f\in\Cal{L}$ such that $\Psi = f\Phi$. In particular,
\begin{equation*}
-\frac{\lambda}{4\gamma} = \opDet(\Psi) = \opDet(\alpha\Phi) = -\frac{\lambda}{4\gamma}f^2.
\end{equation*}
Since $\Cal{L}$ is an algebraic field, it follows that $f =\pm 1$, and uniqueness follows.
\par
We now show that $\Phi$ satisfies \eqnref{PSGuage}. To this end, denote
\begin{equation*}
U := \sum_{k=0}^\infty u_k\lambda^k,\
S := \sum_{k=0}^\infty s_k\lambda^k,\
T := \sum_{k=0}^\infty t_k\lambda^k\ \text{and}\
\Psi := \sum_{k=0}^\infty \psi_k\lambda^k,
\end{equation*}
where $(\psi_m)$ is the sequence constructed in \eqnref{RecursionForPsiI}. By \eqnref{RecursionForTAndS}
\begin{equation*}
T = \frac{1}{\gamma}\big(-2iU_z - \Psi\big)\ \text{and}\
S = \lambda e^{-2\omega}\big(2iU_z - \Psi\big),
\end{equation*}
so that
\begin{equation*}
\eqalign{
Det(\Phi) &= -U^2 - e^{2\omega}ST\vphantom{\frac{1}{2}}
= -U^2+\frac{\lambda}{\gamma}(2iU_z - \Psi)(2iU_z + \Psi)\cr
&= -U^2-4\frac{\lambda}{\gamma} U_z^2 - \frac{\lambda}{\gamma}\Psi^2.\cr}
\end{equation*}
Since the coefficients of the expression
\begin{equation*}
-U^2 -4\frac{\lambda}{\gamma} U_z^2 - \frac{\lambda}{\gamma}\Psi^2 = -\frac{\lambda}{4\gamma}
\end{equation*}
are precisely the recurrence relations \eqnref{RecursionForPsiI} and \eqnref{RecursionForPsiII}, the result follows.\qed
\newsubhead{The Sklyanin matrix}[SklyaninsKMatrix]
We now recall how Sklyanin translates the Durham boundary conditions into boundary conditions for the Lax pair (c.f. \cite{SklyaninI} and \cite{SklyaninII}). Observe first that the real component of the Lax pair is
\begin{equation}
\alpha_x = -\frac{i}{2}\omega_y\sigma_0 + \bigg(\frac{i}{4\lambda}e^\omega + \frac{i}{4\gamma}e^{-\omega}\bigg)\sigma_+ + \bigg(\frac{i\gamma}{4}e^{-\omega} + \frac{i\lambda}{4}e^\omega\bigg)\sigma_-.\eqnum{\nexteqnno[LaxPairReal]}
\end{equation}
We say that $\alpha_x$ satisfies the {\sl Sklyanin condition} for Lax pairs at a point of $\partial\Omega$ whenever
\begin{equation}
K(\lambda,\gamma)\alpha_x(\lambda,\gamma) = \alpha_x\bigg(\frac{1}{\lambda},\frac{1}{\gamma}\bigg)K(\lambda,\gamma)\eqnum{\nexteqnno[KMatrixRelation]}
\end{equation}
at this point, for all $\lambda,\gamma\in\Bbb{C}^*$, where $K(\lambda,\gamma)$ is {\sl Sklyanin's $K$-matrix}, given by \eqnref{KMatrixIntro}. Sklyanin shows
\begin{lemma}\proclabel{SklyaninForLaxHolds}
\noindent The function $\omega$ satisfies the Durham condition
\begin{equation*}
\omega_y = Ae^\omega + Be^{-\omega}
\end{equation*}
at a point of $\partial\Omega$ if and only if the real part $\alpha_x$ of its Lax pair satisfies the Sklyanin condition for Lax pairs at this point.
\end{lemma}
\proof Since the framework of Sklyanin's work is quite different from our own, we include the proof for the reader's convenience. Consider the ansatz
\begin{equation}
K := a\opId + b\sigma_0 + c(\sigma_+ + \sigma_-),\eqnum{\nexteqnno[GeneralK]}
\end{equation}
where the coefficients $a$, $b$ and $c$ only depend on $\gamma$ and $\lambda$. The relation \eqnref{KMatrixRelation} holds if and only if
\begin{equation*}
c\omega_y + \bigg(\frac{(a+b)}{4\lambda} - \frac{(a-b)\lambda}{4}\bigg)e^\omega + \bigg(\frac{(a+b)}{4\gamma} - \frac{(a-b)\gamma}{4}\bigg)e^{-\omega} = 0.
\end{equation*}
Setting
\begin{equation*}
a+b := 4A\gamma - 4B\lambda,\vphantom{\frac{1}{2}},\
a-b := \frac{4A}{\gamma} - \frac{4B}{\lambda}\ \text{and}\
c := \frac{\lambda}{\gamma} - \frac{\gamma}{\lambda},
\end{equation*}
we see that \eqnref{KMatrixRelation} is satisfied at a point of $\partial\Omega$ if and only if
\begin{equation*}
\omega_y = Ae^\omega + Be^{-\omega}
\end{equation*}
at this point, which is precisely the Durham condition. Finally, substituting $a$, $b$ and $c$ into \eqnref{GeneralK} yields \eqnref{KMatrixIntro} as desired.\qed
\newsubhead{The Sklyanin condition}[TheSklyaninCondition]
We conclude this chapter by showing that if $\omega$ satisfies the Durham boundary conditions over $\partial\Omega$, then its Pinkall--Sterling field satisfies the Sklyanin condition for fields over $\partial\Omega$. To this end, let $\partial\Omega_0$ be one of the two connected components of $\partial\Omega$. We define a {\sl Killing field} of $\omega$ over $\partial\Omega_0$ to be a Laurent series $\Phi$ in $\Cal{L}(\partial\Omega_0,\opsl(2,\Bbb{C}))$ which satisfies the Killing field equation
\begin{equation}
\Phi_x = [\Phi,\alpha_x].\eqnum{\nexteqnno[KillingFieldII]}
\end{equation}
Let $\Cal{K}(\partial\Omega_0)$ denote the space of Killing fields of $\omega$ over $\partial\Omega_0$. Trivially, every Killing field of $\omega$ over $\Omega$ restricts to a Killing field of $\omega$ over $\partial\Omega_0$. In the one-dimensional case, we have the following weaker version of Lemma \procref{VanishingConditions}.
\begin{lemma}
\noindent Let
\begin{equation*}
\Phi := \sum_{m=0}^\infty\begin{pmatrix} u_m&e^\omega t_m\cr e^\omega s_m&-u_m\end{pmatrix}\lambda^m
\end{equation*}
be a Killing field of $\omega$ over $\partial\Omega_0$. If
\begin{equation*}
t_{k-2} = s_{k-1} = u_{k-1} = t_{k-1} = 0,
\end{equation*}
then $s_k=u_k=0$ and $t_k$ are constant.
\end{lemma}
\proof By \eqnref{LaxPairReal} and \eqnref{KillingFieldII}, for all $m$,
\begin{align}
u_{m,x} &= \frac{i\gamma}{4}t_m + \frac{i}{4}e^{2\omega}t_{m-1} - \frac{i}{4}e^{2\omega}s_{m+1} - \frac{i}{4\gamma}s_m,\nexteqnno[XRI]\\
e^\omega s_{m,x} &=-2e^\omega\omega_{\overline{z}}s_m - \frac{i\gamma}{2}e^{-\omega}u_m - \frac{i}{2}e^\omega u_{m-1},\ \text{and}\nexteqnno[XRII]\\
e^\omega t_{m,x} &=-2e^\omega\omega_z t_m + \frac{i}{2}e^\omega u_{m+1} + \frac{i}{2\gamma}e^{-\omega}u_m.\nexteqnno[XRIII]
\end{align}
By \eqnref{XRI} and \eqnref{XRIII},
\begin{equation*}
s_k = u_k = 0.
\end{equation*}
These three relations together yield
\begin{equation*}
\eqalign{
e^{2\omega}s_{k+1} &= \gamma t_k,\vphantom{\frac{1}{1}}\cr
u_{k+1} &= \frac{1}{\gamma}\big(2ie^{2\omega}s_{k+1,x}+4ie^{2\omega}\omega_{\overline{z}}s_{k+1}\big),\ \text{and}\cr
e^\omega t_{k,x} &= -2e^\omega\omega_z t_k + \frac{i}{2}e^\omega u_{k+1}.\vphantom{\frac{1}{2}}\cr}
\end{equation*}
It follows that
\begin{equation*}
u_{k+1} = 2ie^{2\omega}(e^{-2\omega}t_k)_{,x} + 4i\omega_{\overline{z}}t_k=-4i\omega_xt_k + 2it_{k,x} + 4i\omega_{\overline{z}}t_k,
\end{equation*}
so that
\begin{equation*}
e^\omega t_{k,x} = -2e^\omega\omega_z t_k + 2e^\omega\omega_x t_k - e^\omega t_{k,x} - 2e^\omega\omega_{\overline{z}}t_k.
\end{equation*}
It follows that $t_{k,x}=0$, and this completes the proof.\qed
\\ \null \\
\noindent Via the same argument as in Sections \subheadref{LEspaceDeChampsDeKilling} and \subheadref{LeDeterminant}, this yields
\begin{lemma}\proclabel{KillingFieldsAlongBoundary}
\myitem{(1)} $\Cal{K}(\partial\Omega_0)$ is the one-dimensional vector space over $\Cal{L}$ generated by the Pinkall--Sterling field, and
\medskip
\myitem{(2)} the restriction of the Pinkall--Sterling field to $\partial\Omega_0$ is, up to a choice of sign, the unique Killing field $\Phi$ of $\omega$ over $\partial\Omega_0$ which satisfies
\begin{equation*}
\opDet(\Phi) = -\frac{\lambda}{4\gamma}.
\end{equation*}
\end{lemma}
Given a Killing field $\Phi$ of $\omega$ over $\partial\Omega_0$, denote
\begin{equation}
\tilde{\Phi}(\lambda,\gamma) := K(\lambda,\gamma)^{-1}\overline{\Phi(\overline{\lambda},\overline{\gamma})}^tK(\lambda,\gamma).\eqnum{\nexteqnno[PhiTilde]}
\end{equation}
Observe that $\tilde{\Phi}$ is also a Laurent series over $\partial\Omega_0$.
\begin{lemma}\proclabel{ConjugationByKPreservesKilling}
\noindent If the real part $\alpha_x$ of the Lax pair of $\omega$ satisfies the Sklyanin condition for Lax pairs along $\partial\Omega_0$, then $\tilde{\Phi}$ satisfies the Killing field equation \eqnref{KillingFieldII} along $\partial\Omega_0$.
\end{lemma}
\proof Observe that
\begin{equation*}
\eqalign{
\overline{\alpha_x(\overline{\lambda}^{-1},\overline{\gamma}^{-1})}^t &= -\alpha_x(\lambda,\gamma),\vphantom{\bigg|}\cr
\overline{K(\overline{\lambda}^{-1},\overline{\gamma}^{-1})} &= D(\lambda,\gamma)K(\lambda,\gamma)^{-1}\ \text{and}\vphantom{\bigg|}\cr
K(\lambda^{-1},\gamma^{-1}) &= D(\lambda,\gamma)K(\lambda,\gamma)^{-1},\cr}
\end{equation*}
where
\begin{equation*}
D(\lambda,\gamma) := \opDet(K(\lambda,\gamma)).
\end{equation*}
Moreover,
\begin{equation*}
D(\lambda,\gamma) = D(\lambda^{-1},\gamma^{-1}) = \overline{D(\overline{\lambda}^{-1},\overline{\gamma}^{-1})}.
\end{equation*}
Upon applying the Killing field equation \eqnref{KillingFieldII}, we therefore obtain
\begin{equation*}
\eqalign{
\tilde{\Phi}(\lambda,\gamma)_x &= K(\lambda,\gamma)^{-1}\overline{\Phi(\overline{\lambda},\overline{\gamma})_x}^tK(\lambda,\gamma)\cr
&= K(\lambda,\gamma)^{-1}\overline{\big[\Phi(\overline{\lambda},\overline{\gamma}),\alpha_x(\overline{\lambda},\overline{\gamma})\big]}^tK(\lambda,\gamma)\cr
&= K(\lambda,\gamma)^{-1}\big[\overline{\alpha_x(\overline{\lambda},\overline{\gamma})}^t,\overline{\Phi(\overline{\lambda},\overline{\gamma})}^t\big]K(\lambda,\gamma)\cr
&= -K(\lambda,\gamma)^{-1}\big[\alpha_x(\lambda^{-1},\gamma^{-1}),\overline{\Phi(\overline{\lambda},\overline{\gamma})}^t\big]K(\lambda,\gamma)\cr
&= K(\lambda,\gamma)^{-1}\big[\overline{\Phi(\overline{\lambda},\overline{\gamma})}^t,\alpha_x(\lambda^{-1},\gamma^{-1})\big]K(\lambda,\gamma)\cr
&= \big[K(\lambda,\gamma)^{-1}\overline{\Phi(\overline{\lambda},\overline{\gamma})}^tK(\lambda,\gamma),K(\lambda,\gamma)^{-1}\alpha_x(\lambda^{-1},\gamma^{-1})K(\lambda,\gamma)\big]\cr
&= [\tilde{\Phi}(\lambda,\gamma),\alpha_x(\lambda,\gamma)],\vphantom{\big[}\cr}
\end{equation*}
and the result follows.\qed
\begin{lemma}\proclabel{PinkalSterlingFieldIsReal}
\noindent If the real part $\alpha_x$ of the Lax pair of $\omega$ satisfies the Sklyanin condition for Lax pairs along $\partial\Omega_0$, then the Pinkall--Sterling field $\Phi$ of $\omega$ satisfies the Sklyanin condition for fields along $\partial\Omega_0$, that is, for all $\lambda,\gamma\in S^1$,
\begin{equation}
\Phi(\lambda,\gamma) = \tilde{\Phi}(\lambda,\gamma) = K(\lambda,\gamma)^{-1}\overline{\Phi(\overline{\lambda},\overline{\gamma})}^tK(\lambda,\gamma)\eqnum{\nexteqnno[KReality]}
\end{equation}
along $\partial\Omega_0$.
\end{lemma}
\proof By Lemma \procref{ConjugationByKPreservesKilling}, $\tilde{\Phi}$ is a Killing field of $\omega$ over $\partial\Omega_0$. However, for all $\lambda$ and for all $\gamma$,
\begin{equation*}
\opDet(\tilde{\Phi}(\lambda,\gamma))
=\opDet\big(K(\lambda,\gamma)^{-1}\overline{\Phi(\overline{\lambda},\overline{\gamma})}^tK(\lambda,\gamma)\big)
=\overline{\opDet\big(\Phi(\overline{\lambda},\overline{\gamma})\big)}
=-\frac{\lambda}{4\gamma}.
\end{equation*}
It follows by Lemma \procref{KillingFieldsAlongBoundary} that
\begin{equation*}
\Phi = \pm\tilde{\Phi}.
\end{equation*}
The result now follows upon explicitly calculating the first non-zero term of each of these two series.\qed
%
%\newhead{The finite-type property}[TheFiniteTypeProperty]
%
\newsubhead{Robin boundary conditions}[LaConditionDeBordDeRobin]
As before, let $\Omega:=\Bbb{R}\times[-T,T]$ and let $\omega:\Omega\rightarrow\Bbb{R}$ be a real solution of the elliptic sinh--Gordon equation. In this chapter, we transform the Sklyanin condition for fields into a sequence of equations that allow us to recover the finite type property. Thus, let $\Phi$ be the Pinkall--Sterling field of $\omega$. In order to better capture the symmetries of the problem, it now becomes convenient to treat this field as a Laurent series in $\lambda$ and $\gamma$. We thus denote
\begin{equation}
\Phi(\lambda,\gamma) := \sum_{m,n}\begin{pmatrix} u_{m,n}& e^\omega t_{m,n}\cr e^\omega s_{m,n} & -u_{m,n}\cr\end{pmatrix}\lambda^m\gamma^n.\eqnum{\nexteqnno[NewExpansionForPhi]}
\end{equation}
Observe that, for all $\lambda,\gamma\in S^1$,
\begin{equation}
\alpha(\lambda,\gamma) = e^{-\frac{\theta}{2}\sigma_0}\alpha\bigg(\frac{\lambda}{\gamma},1\bigg)e^{\frac{\theta}{2}\sigma_0},\eqnum{\nexteqnno[GuageTransform]}
\end{equation}
where $\theta\in\Bbb{R}$ satisfies
\begin{equation*}
e^{2i\theta}=\gamma.
\end{equation*}
It follows upon applying this gauge transformation that
\begin{equation}
\Phi(\lambda,\gamma) = \sum_{m=0}^\infty\begin{pmatrix} u_m&\gamma^{-1}e^\omega t_m\cr \gamma e^\omega s_m& -u_m\cr\end{pmatrix}\lambda^m\gamma^{-m},\eqnum{\nexteqnno[PSFieldGeneralGamma]}
\end{equation}
where
\begin{equation*}
\sum_{m=0}^\infty\begin{pmatrix} u_m& e^\omega t_m\cr e^\omega s_m& -u_m\cr\end{pmatrix}\lambda^m := \Phi(\lambda,1)
\end{equation*}
is the Pinkall--Sterling field of $\omega$ with torsion $\gamma=1$. In this framework, Lemma \procref{LaxFormulae} becomes
\begin{lemma}
\noindent The sequences $(u_{m,n})$, $(t_{m,n})$ and $(s_{m,n})$ satisfy, for all $m$ and for all $n$,
\begin{align}
4u_{m,n,z} + i e^{2\omega} s_{m+1,n} - i t_{m,n-1} &=0, \nexteqnno[NRI]\\
4u_{m,n,\overline{z}} + is_{m,n+1} - ie^{2\omega} t_{m-1,n} &=0, \nexteqnno[NRII]\\
4\omega_z t_{m,n} + 2 t_{m,n,z} - i u_{m+1,n} &=0, \nexteqnno[NRIII]\\
2e^\omega t_{m,n,\overline{z}} - i e^{-\omega}u_{m,n+1} &=0, \nexteqnno[NRIV]\\
2e^\omega s_{m,n,z} + i e^{-\omega} u_{m,n-1} &=0\ \text{and} \nexteqnno[NRV]\\
4\omega_{\overline{z}}s_{m,n} + 2s_{m,n,\overline{z}} + iu_{m-1,n} &=0.\nexteqnno[NRVI]
\end{align}
\end{lemma}
\noindent As far as the Sklyanin condition for fields is concerned, upon equating every coefficient of $(K\Phi - \tilde{\Phi}K)$ with zero, we obtain
\begin{lemma}
\noindent Along $\partial\Omega$, the sequences $(u_{m,n})$, $(t_{m,n})$ and $(s_{m,n})$ satisfy, for all $m$ and for all $n$,
\begin{align}
\opIm\big(e^\omega s_{m-1,n+1} - e^\omega s_{m+1,n-1} + 4Au_{m,n-1} - 4Bu_{m-1,n}\big) &= 0,\nexteqnno[SI]\\
\opIm\big(e^\omega t_{m-1,n+1} - e^\omega t_{m+1,n-1} - 4Au_{m,n+1} + 4Bu_{m+1,n}\big) &= 0,\nexteqnno[SII]\\
\opRe\big(2Ae^\omega t_{m,n-1} - 2Be^\omega t_{m-1,n} - 2Ae^\omega s_{m,n+1} + 2Be^\omega s_{m+1,n}\big) & \nexteqnno[SIII]\\
\qquad\qquad\qquad=\opRe\big(u_{m-1,n+1} - u_{m+1,n-1}\big)\ \text{and}& \notag\\
\opIm\big(At_{m,n-1} - Bt_{m-1,n} + As_{m,n+1} - Bs_{m+1,n}\big) &= 0.\nexteqnno[SIV]
\end{align}
\end{lemma}
\begin{lemma}
\noindent The pair of systems of equations \eqnref{SI} and \eqnref{SII} is equivalent to the following pair of systems of equations:
\begin{align}
\opIm\big(e^\omega t_{m-1,n} - 4Au_{m,n} + e^\omega s_{m+1,n}\big) &= 0\ \text{and}\nexteqnno[NSI]\\
\opIm\big(e^\omega t_{m,n-1} - 4Bu_{m,n} + e^\omega s_{m,n+1}\big) &= 0.\nexteqnno[NSII]
\end{align}
\end{lemma}
\proof Indeed, suppose that the equations \eqnref{SI} and \eqnref{SII} are satisfied for all $m$ and for all $n$. Then, upon applying recursively \eqnref{SI}, we obtain the {\sl finite} sum
\begin{equation*}
\opIm\big(e^\omega s_{m,n}\big) = \opIm\big(4A u_{m-1,n} - 4B u_{m-2,n+1} + 4Au_{m-3,n+2} - \cdots\big).
\end{equation*}
In the same manner, \eqnref{SII} yields
\begin{equation*}
\opIm\big(e^\omega t_{m,n}\big) = \opIm\big(4B u_{m,n+1} - 4A u_{m-1,n+2} + 4B u_{m-2,n+3} - \cdots\big).
\end{equation*}
If we now denote these equations respectively by $\alpha(m,n)$ and $\beta(m,n)$, then $\alpha(m+1,n)$ and $\beta(m-1,n)$ together yield \eqnref{NSI} whilst $\alpha(m,n+1)$ and $\beta(m,n-1)$ together yield \eqnref{NSII}. Since the converse is trivial, this completes the proof.\qed
\begin{lemma}\proclabel{BoundaryConditionsForImaginaryPart}
\noindent For all $m$ and for all $n$, the imaginary part of the function $u_{m,n}$ satisfies the following Robin boundary condition:
\begin{equation}
\opIm(u_{m,n})_y = Ae^\omega\opIm(u_{m,n}) - Be^{-\omega}\opIm(u_{m,n}).\eqnum{\nexteqnno[RobinBoundaryCondition]}
\end{equation}
\end{lemma}
\proof Indeed, by \eqnref{NRI} and \eqnref{NRII},
\begin{equation*}
\eqalign{
e^\omega s_{m+1,n} - e^{-\omega}t_{m,n-1} &= 4ie^{-\omega}u_{m,n,z}\ \text{and}\cr
e^\omega t_{m-1,n} - e^{-\omega}s_{m,n+1} &= -4ie^{-\omega}u_{m,n,\overline{z}}.\cr}
\end{equation*}
The sum of \eqnref{NSI} and \eqnref{NSII} then yields
\begin{equation*}
\eqalign{
\opIm\big(4A u_{m,n} - 4B e^{-2\omega} u_{m,n}\big)
&=\opIm\big(e^\omega t_{m-1,n} + e^\omega s_{m+1,n} - e^{-\omega} t_{m,n-1} - e^{-\omega} s_{m,n+1}\big)\cr
&=\opIm\big(4ie^{-\omega}u_{m,n,z} - 4ie^{-\omega}u_{m,n,\overline{z}}\big)\cr
&=4e^{-\omega}\opRe\big((u_{m,n} - \overline{u}_{m,n})_z\big)\cr
&=4e^{-\omega}\opRe\big(2i\opIm(u_{m,n})_z\big)\cr
&=4e^{-\omega}\opIm(u_{m,n})_y,\cr}
\end{equation*}
and the result follows.\qed
\newsubhead{The solution is of finite type}[LaSolutionEstDeTypeFini]
We now show that the solution is of finite type. Let $\Phi$ be as in the previous section. Recall that, upon applying the guage transformation \eqnref{GuageTransform} if necessary, we may henceforth suppose that $\gamma=1$. We first recall a few elementary lemmas.
\begin{lemma}\proclabel{PreFiniteType}
\noindent Let
\begin{equation*}
\Psi := \sum_{m=0}^\infty\begin{pmatrix} u_m&e^\omega t_m\cr e^\omega s_m&-u_m\cr\end{pmatrix}\lambda^m
\end{equation*}
be a Killing field. If $u_k=0$ then there exists a Laurent series $f\in\Cal{L}$ of degree $k$ such that
\begin{equation*}
\Psi - f\Phi = \sum_{m=0}^{k-1}\begin{pmatrix} u_m&e^\omega t_m\cr e^\omega s_m&-u_m\cr\end{pmatrix}\lambda^m + \begin{pmatrix} 0&0\cr e^\omega s_k&0\cr\end{pmatrix}\lambda^k.
\end{equation*}
In particular, $(\Psi-f\Phi)$ is a polynomial Killing field.
\end{lemma}
\proof We construct $f$ by recurrence. Suppose that the coefficients $f_k,\cdots,f_{k+l-1}$ have already been determined such that if
\begin{equation*}
\tilde{\Psi} := \Psi - f_{(l)}\Phi := \sum_{m=0}^\infty\begin{pmatrix}\tilde{u}_{m}&e^\omega\tilde{t}_{m}\cr e^\omega\tilde{s}_{m}&-\tilde{u}_{m}\end{pmatrix}\lambda^m,
\end{equation*}
where
\begin{equation*}
f_{(l)} := \sum_{m=k}^{k+l-1}f_m\lambda^m,
\end{equation*}
then,
\begin{equation*}
\eqalign{
\tilde{u}_{m} &=0\ \forall k\leq m\leq k+l,\cr
\tilde{t}_{m} &=0\ \forall k\leq m\leq k+l-1\ \text{and}\cr
\tilde{s}_{m} &=0\ \forall k+1\leq m\leq k+l.\cr}
\end{equation*}
Since $\tilde{\Psi}$ is also a Killing field, by Lemma \procref{VanishingConditions},
\begin{equation*}
\tilde{t}_{k+l} = c
\end{equation*}
is constant. By Lemma \procref{VanishingConditions} again, the result follows upon setting $f_{k+l}:=-2c$.\qed
\begin{lemma}\proclabel{FiniteTypeCriteria}
\noindent The solution $\omega$ is of finite type if and only if there exists a finite-dimensional vector space $E\subseteq C^\infty(\Omega,\Bbb{C})$ such that, for all $m$ and for all $n$,
\begin{equation*}
u_{m,n}\in E.
\end{equation*}
\end{lemma}
\proof This condition is trivially necessary. We now show that it is sufficient. Suppose again that $\gamma=1$. In particular, by \eqnref{PSFieldGeneralGamma}, for all $m$, $u_m = u_{m,m} \in E$. If $d:=\opDim(E)$, then there exists a Laurent polynomial $f\in\Cal{P}$ of bidegree $(0,d-1)$ such that if
\begin{equation*}
f\Phi := \sum_{m=0}^\infty\begin{pmatrix}\tilde{u}_m& e^\omega\tilde{t}_m\cr e^\omega\tilde{s}_m&-\tilde{u}_m\end{pmatrix}\lambda^m,
\end{equation*}
then $\tilde{u}_d = 0$, and the result follows by Lemma \procref{PreFiniteType}.\qed
\begin{lemma}
\noindent For all $m$ and for all $n$, the function $u_{m,n}$ satisfies the following linearised sinh--Gordon equation:
\begin{equation}
\Delta u_{m,n} + \opCosh(2\omega)u_{m,n} = 0.\eqnum{\nexteqnno[LinearisedSinhGordon]}
\end{equation}
\end{lemma}
\proof Indeed, differentiating \eqnref{NRI} yields
\begin{equation*}
4u_{m,n,z\overline{z}} + i(e^{2\omega}s_{m+1,n})_{\overline{z}} - it_{m,n-1,\overline{z}} = 0.
\end{equation*}
Applying \eqnref{NRIV} and \eqnref{NRVI} then yields
\begin{equation*}
4u_{m,n,z\overline{z}} + \frac{1}{2}e^{2\omega}u_{m,n} + \frac{1}{2}e^{-2\omega} u_{m,n} = 0,
\end{equation*}
as desired.\qed
\begin{lemma}\proclabel{Liouville}
\noindent If $\phi:\Omega\rightarrow\Bbb{C}$ is a periodic, holomorphic function which satisfies
\begin{equation*}
\opIm(\phi)|_{\partial\Omega} = 0,
\end{equation*}
then $\phi$ is constant.
\end{lemma}
\proof Indeed, by Cauchy's reflection principle, $\phi$ extends to a bounded holomorphic function over $\Bbb{C}$ and the result now follows by Liouville's theorem.\qed
\begin{theorem}
\noindent If $\omega:\Omega\rightarrow\Bbb{R}$ is a periodic solution of the sinh--Gordon equation with Durham boundary conditions, then $\omega$ is of finite type.
\end{theorem}
\proof Indeed, let $\Phi$ be the Pinkall--Sterling field of $\omega$ and let $(u_{m,n})$, $(t_{m,n})$ and $(s_{m,n})$ be as in \eqnref{NewExpansionForPhi}. By Lemmas \procref{SklyaninForLaxHolds} and \procref{PinkalSterlingFieldIsReal}, $\Phi$ satisfies the Sklyanin condition for fields along $\partial\Omega$. It follows by Lemma \procref{BoundaryConditionsForImaginaryPart} that, for all $m$ and for all $n$,
\begin{equation*}
\opIm(u_{m,n})_{,y} = Ae^\omega\opIm(u_{m,n}) - Be^{-\omega}\opIm(u_{m,n}),
\end{equation*}
along $\partial\Omega$. Since $\opIm(u_{m,n})$ also satisfies the linearised sinh--Gordon equation, it follows by the classical theory of elliptic operators over compact manifolds with boundary (see \cite{GilbTrud}) that there exists a finite-dimensional subspace $E_1\subseteq C^\infty(\Omega,\Bbb{R})$ such that, for all $m$ and for all $n$,
\begin{equation*}
\opIm(u_{m,n}) \in E_1.
\end{equation*}
\par
By \eqnref{NRIII} and \eqnref{NRVI}, for all $m$ and for all $n$,
\begin{equation*}
\big(e^{2\omega}\overline{s}_{m+1,n} - e^{2\omega}t_{m-1,n}\big)_z
=-\frac{i}{2}e^{2\omega}\big(u_{m,n} - \overline{u}_{m,n}\big)=e^{2\omega}\opIm(u_{m,n})\in e^{2\omega}E_1,
\end{equation*}
and, by \eqnref{NSI}, along $\partial\Omega$,
\begin{equation*}
\opIm\big(e^{2\omega}\overline{s}_{m+1,n} - e^{2\omega}t_{m-1,n}\big)=-4Ae^\omega\opIm(u_{m,n})\in e^\omega E_1|_{\partial\Omega}.
\end{equation*}
It follows by Lemma \procref{Liouville} that there exists a finite-dimensional subspace $E_2\subseteq C^\infty(\Omega,\Bbb{C})$ such that, for all $m$ and for all $n$,
\begin{equation*}
e^{2\omega}\overline{s}_{m+1,n} - e^{2\omega}t_{m-1,n} \in e^\omega E_2.
\end{equation*}
In a similar manner, we show by \eqnref{NRIV} and \eqnref{NRV}, that there exists a finite-dimensional subspace $E_3\subseteq C^\infty(\Omega,\Bbb{C})$ such that, for all $m$ and for all $n$,
\begin{equation*}
s_{m,n+1} - \overline{t}_{m,n-1} \in e^{-\omega}E_3.
\end{equation*}
\par
Finally, by \eqnref{SIII}, along $\partial\Omega$, for all $m$ and for all $n$,
\begin{equation*}
\eqalign{
\opRe\big(u_{m+1,n-1} - u_{m-1,n+1}\big)
&=2Be^{-\omega}\opRe\big(e^{2\omega}\overline{s}_{m+1,n}-e^{2\omega}t_{m-1,n}\big)\cr
&\qquad\quad-2Ae^\omega\opRe\big(s_{m,n+1} - \overline{t}_{m,n-1}\big)\cr
&\in\opRe(E_2+E_3).\cr}
\end{equation*}
It follows by induction that, along $\partial\Omega$, for all $m$ and for all $n$,
\begin{equation*}
\opRe(u_{m,n}) \in \opRe(E_2+E_3).
\end{equation*}
Finally, since $\opRe(u_{m,n})$ also satisfies the linearised sinh--Gordon equation, it follows again by the classical theory of elliptic operators over compact manifolds with boundary that there exists a fourth finite-dimensional subspace $E_4\subseteq C^\infty(\Omega,\Bbb{R})$ such that, for all $m$ and for all $n$,
\begin{equation*}
\opRe(u_{m,n}) \in E_4.
\end{equation*}
The result now follows by Lemma \procref{FiniteTypeCriteria}.\qed
\newsubhead{Polynomial Killing fields}[DesChampsDeKillingPolynomes]
Finally, we construct polynomial Killing fields of $\omega$ over $\Omega$ that also satisfy the Sklyanin condition for fields. As in Section \subheadref{LaSolutionEstDeTypeFini} restrict ourselves again to the case where $\gamma=1$.
\begin{lemma}\proclabel{ExceptionalKillingField}
\noindent Let
\begin{equation*}
\Psi := \sum_{m=0}^\infty\begin{pmatrix} u_m&e^\omega t_m\cr e^\omega s_m&- u_m\cr\end{pmatrix}\lambda^m
\end{equation*}
be a Killing field of $\omega$ over $\Omega$ which satisfies the Sklyanin condition for fields. If $u_k=0$ and if $t_k\notin\Bbb{R}$ then there exists a polynomial Killing field of $\omega$ over $\Omega$ of bidegree $(0,4)$ which also satisfies the Sklyanin condition for fields.
\end{lemma}
\proof Let $\Phi$ be the Pinkall--Sterling field. By Lemma \procref{PreFiniteType}, there exist Laurent series $f,g\in\Cal{L}$ with real coefficients such that
\begin{equation*}
\Psi - (f+ig)\Phi = P,
\end{equation*}
where $P$ is a polynomial Killing field of $\omega$ over $\Omega$. We therefore denote
\begin{equation*}
\Psi_1 := g\Phi.
\end{equation*}
Since $\Phi$ satisfies the Sklyanin condition for fields, and since $g$ has real coefficients, $\Psi_1$ also satisfies the Sklyanin condition for fields, that is, for all $\lambda\in\Bbb{C}^*$,
\begin{equation*}
K(\lambda,1)\Psi_1(\lambda,1) - \overline{\Psi_1(\overline{\lambda},1)}^tK(\lambda,1)=0
\end{equation*}
along $\partial\Omega$. Next, since $\Phi$ and $\Psi$ both satisfy the Sklyanin condition for fields, and since $f$ has real coefficients, we also have, for all $\lambda\in S^1$,
\begin{equation*}
iK(\lambda,1)\Psi_1(\lambda,1)+i\overline{\Psi_1(\overline{\lambda},1)}^tK(\lambda,1)=Q(\lambda,1)
\end{equation*}
along $\partial\Omega$, where
\begin{equation*}
Q(\lambda,1):=K(\lambda,1)P(\lambda,1) - \overline{P(\overline{\lambda},1)}^tK(\lambda,1).
\end{equation*}
It follows that
\begin{equation*}
\triplealign{
&K(\lambda,1)\Psi_1(\lambda,1) &= -\frac{i}{2}Q(\lambda,1)\cr
\Leftrightarrow &D(\lambda,1)\Psi_1(\lambda,1) &= -\frac{i}{2}K(\lambda^{-1},1)Q(\lambda,1),\cr}
\end{equation*}
where
\begin{equation*}
D(\lambda,1) := \opDet(K(\lambda,1)).
\end{equation*}
We thus denote
\begin{equation*}
\Psi'(\lambda,1) := \lambda^{2-k}D(\lambda,1)\Psi_1(\lambda,1).
\end{equation*}
Since $D(\lambda,1)$ is a Laurent polynomial with real coefficients, $\Psi'$ is also a Killing field of $\omega$ over $\Omega$ which satisfies the Sklyanin condition for fields. Finally, we verify that $Q(\lambda,1)$ is a Laurent polynomial of bidegree $(k-1,k+1)$, and the result follows.\qed
\\ \null \\
{\bf\noindent Proof of Theorem \procref{MainTheorem}:}\ In order to prove the existence of a polynomial Killing field of $\omega$ over $\Omega$ which satisfies the Sklyanin condition for fields, it suffices to repeat the construction of Lemmas \procref{PreFiniteType} and \procref{FiniteTypeCriteria} using only Laurent series with real coefficients. The only possible obstruction to this construction is precisely the case studied in Lemma \procref{ExceptionalKillingField}, where there nonetheless exists a polynomial Killing field of $\omega$ over $\Omega$ which satisfies the Sklyanin condition for fields. This proves existence.
\par
Conversely, suppose that there exists a polynomial Killing field of $\omega$ over $\Omega$
\begin{equation*}
\Psi := \sum_{m=0}^k\begin{pmatrix} u_m&e^\omega t_m\cr e^\omega s_m&-u_m\cr\end{pmatrix}\lambda^m
\end{equation*}
which satisfies the Sklyanin condition for fields. Then, by Lemma \procref{VanishingConditions}, $s_0 = u_0 = 0$. By \eqnref{SII}, we may suppose that $t_0 = 1/2$ and, by \eqnref{SIII},
\begin{equation*}
\opRe\big(u_1 + 2Ae^\omega t_0 + 2Be^\omega s_1\big) = 0.
\end{equation*}
However, by \eqnref{InitialValuesOfRecursion}, \eqnref{RecursionForU} and \eqnref{RecursionForTAndS},
\begin{equation*}
u_1 = -2i\omega_z,\vphantom{\frac{1}{2}},\
t_0 =\frac{1}{2}\ \text{and}\
s_0 =\frac{1}{2}e^{-2\omega}.
\end{equation*}
The result follows.\qed
\bibliographystyle{amsplain}
\def\cydot{\leavevmode\raise.4ex\hbox{.}} \def\cprime{$'$}
\providecommand{\bysame}{\leavevmode\hbox to3em{\hrulefill}\thinspace}
\providecommand{\MR}{\relax\ifhmode\unskip\space\fi MR }
\providecommand{\MRhref}[2]{%
  \href{http://www.ams.org/mathscinet-getitem?mr=#1}{#2}
}
\providecommand{\href}[2]{#2}

\end{document}